\documentclass{article} 
\usepackage{mathmag_arxiv}  

\usepackage{tikz} 
\usetikzlibrary{positioning,decorations.pathmorphing} 
\usepackage{booktabs}
\usepackage{xcolor}

\usepackage{graphicx}
\usepackage{chngcntr}

\theoremstyle{theorem}

\theoremstyle{definition}

\allowdisplaybreaks

\makeatletter
\@addtoreset{footnote}{page}
\makeatother

\begin{document}

\title{Counting Connected and Disconnected Ways to Assemble a Jigsaw Puzzle}

\author{
Prarthana Agrawal\\               
\scriptsize Rudolf Peierls Centre for Theoretical Physics\\University of Oxford\\ 
prarthana.agrawal@physics.ox.ac.uk\\                   
\medskip                                          
\normalsize{Abdurrahman Hadi Erturk}\\         
\scriptsize Rudolf Peierls Centre for Theoretical Physics\\University of Oxford\\    
abdurrahman.erturk@physics.ox.ac.uk\\
\medskip                                         
\normalsize{Ard A. Louis}\\        
\scriptsize Rudolf Peierls Centre for Theoretical Physics\\University of Oxford\\ 
ard.louis@physics.ox.ac.uk
}                               

\maketitle

\begin{abstract}
A jigsaw puzzle may be assembled in many different ways. Some assembly sequences remain connected throughout, while others temporarily build separate parts of the puzzle before joining them together. By representing the puzzle as a graph, these assembly sequences become vertex orderings that can be counted using recent results from graph theory. We apply this framework to enumerate three natural assembly strategies: connected assembly, assembly beginning from several disconnected pieces, and assembly in which new disconnected sections may be started during the process. The resulting counts reveal that, even for modest puzzle sizes, connectivity-preserving assembly sequences are greatly outnumbered by those that pass through disconnected intermediate stages.
\end{abstract}


\section{Introduction}
\label{sec:introduction}

\begin{quote}
How many different ways are there to assemble a jigsaw puzzle?
\end{quote}

\begin{figure}[t]
\centering
\noindent
\begin{minipage}[t]{0.32\textwidth}
\centering
    \begin{tikzpicture}[scale=0.60, every node/.style={font=\small}]
    \begin{scope}
    \clip (0,0) rectangle (5,5);
    \node[anchor=center, inner sep=0] at (2.55,0.8)
    {
    \includegraphics[
    width=7cm,
    height=5cm,
    keepaspectratio=true,
    trim=0 0 0 0,
    clip
    ]{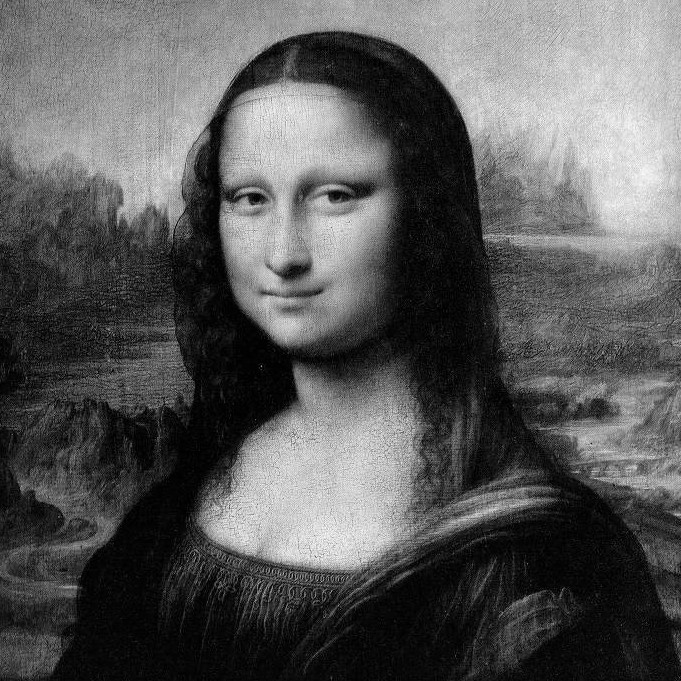}
    };
    \end{scope}
    \foreach \i in {0,...,4} {
    \foreach \j in {0,...,4} {
    \draw[very thick, rounded corners=1pt] (\i,\j) rectangle +(1,1);
    }
    }
    \draw[step=1cm, white!50, thin] (0,0) grid (5,5);
    \end{tikzpicture}
\vspace{1mm}

{\small (a) $5\times 5$ Jigsaw puzzle}
\end{minipage}
\hfill
\begin{minipage}[t]{0.32\textwidth}
\centering
    \begin{tikzpicture}[scale=0.60, every node/.style={font=\small}]
    
    \begin{scope}
    \clip (0,0) rectangle (5,5);
    \node[anchor=center, inner sep=0, opacity=0.18] at (2.55,0.8)
    {
    \includegraphics[
    width=7cm,
    height=5cm,
    keepaspectratio=true,
    trim=0 0 0 0,
    clip
    ]{mona_lisa.jpeg}
    };
    \end{scope}
    \foreach \x/\y in {1/1,2/1,2/2,3/2,3/3} {
    \begin{scope}
    \clip (\x,\y) rectangle +(1,1);
    \node[anchor=center, inner sep=0, opacity=0.90] at (2.55,0.8)
    {
    \includegraphics[
    width=7cm,
    height=5cm,
    keepaspectratio=true,
    trim=0 0 0 0,
    clip
    ]{mona_lisa.jpeg}
    };
    \end{scope}
    }
    \foreach \i in {0,...,4} {
    \foreach \j in {0,...,4} {
    \draw[very thick, rounded corners=1pt] (\i,\j) rectangle +(1,1);
    }
    }
    \draw[step=1cm, white!45, thin] (0,0) grid (5,5);
    \node[text=white] at (2.5,1.5) {\small \textbf{1}};
    \node[text=white] at (1.5,1.5) {\small \textbf{3}};
    \node[text=white] at (2.5,2.5) {\small \textbf{2}};
    \node[text=white] at (3.5,2.5) {\small \textbf{4}};
    \node[text=white] at (3.5,3.5) {\small \textbf{5}};
    \end{tikzpicture}
\vspace{1mm}

{\small (b) Connected assembly}
\end{minipage}
\hfill
\begin{minipage}[t]{0.32\textwidth}
    \centering
    \begin{tikzpicture}[scale=0.60, every node/.style={font=\small}]
    \begin{scope}
    \clip (0,0) rectangle (5,5);
    \node[anchor=center, inner sep=0, opacity=0.18] at (2.55,0.8)
    {
    \includegraphics[
    width=7cm,
    height=5cm,
    keepaspectratio=true,
    trim=0 0 0 0,
    clip
    ]{mona_lisa.jpeg}
    };
    \end{scope}
    \foreach \x/\y in {0/0,1/0,1/1,3/3,3/4,4/3} {
    \begin{scope}
    \clip (\x,\y) rectangle +(1,1);
    \node[anchor=center, inner sep=0, opacity=0.90] at (2.55,0.8)
    {
    \includegraphics[
    width=7cm,
    height=5cm,
    keepaspectratio=true,
    trim=0 0 0 0,
    clip
    ]{mona_lisa.jpeg}
    };
    \end{scope}
    }
    \foreach \i in {0,...,4} {
    \foreach \j in {0,...,4} {
    \draw[very thick, rounded corners=1pt] (\i,\j) rectangle +(1,1);
    }
    }
    \draw[step=1cm, white!45, thin] (0,0) grid (5,5);
    \node[text=white] at (0.5,0.5) {\small \textbf{1}};
    \node[text=white] at (1.5,0.5) {\small \textbf{2}};
    \node[text=white] at (1.5,1.5) {\small \textbf{4}};
    \node[text=white] at (3.5,3.5) {\small \textbf{3}};
    \node[text=white] at (3.5,4.5) {\small \textbf{5}};
    \node[text=white] at (4.5,3.5) {\small \textbf{6}};
    \end{tikzpicture}
\vspace{1mm}

{\small (c) Disconnected assembly}
\end{minipage}
\caption{
Illustration of connected and disconnected assembly on a $5\times5$ jigsaw puzzle. (a) The completed puzzle. (b) A connected assembly sequence, in which every newly placed piece is adjacent to the existing partial assembly. (c) A disconnected assembly sequence, in which the third piece creates a new disconnected component. The numbers indicate the order in which the illustrated pieces are placed.
}
\label{fig:mona_lisa}
\end{figure}

At first glance the answer seems obvious: if the puzzle has $n$ pieces, there are $n!$ possible orders in which they could be placed, a number that grows very fast with $n$. For example, the $5\times 5$ puzzle in Figure~\ref{fig:mona_lisa}a has $25! \approx 1.55 \times 10^{25}$ different ways of solving it. Most of these $n!$ ways, however, bear little resemblance to how a puzzle is actually assembled. A natural strategy is to begin from a single piece and repeatedly attach new pieces so that the partially assembled puzzle remains connected throughout (Figure~\ref{fig:mona_lisa}b). A more relevant question, then, is: how many assembly sequences remain connected throughout the assembly process?

Connected assembly, however, is only one possibility. Many solvers instead begin from several distant pieces, such as the corners, and grow each partial assembly independently before eventually joining them together. An even more general strategy is to place distant pieces whenever they are identified, allowing new disconnected islands
to appear at arbitrary stages during the assembly process (Figure~\ref{fig:mona_lisa}c). Assembly then continues by growing all disconnected islands until they eventually merge into the completed puzzle. How many assembly sequences are possible under these more general assembly strategies?

In this article, we show that this problem can be solved by representing a jigsaw puzzle as a graph, in which each puzzle piece corresponds to a vertex and two vertices are joined whenever the corresponding pieces share a border. An assembly sequence then becomes a linear ordering of the vertices. Requiring the partially assembled puzzle to remain connected after every step is equivalent to requiring that the vertices placed so far in the ordering induce a connected subgraph. In graph-theoretic language, such orderings are known as \emph{successive vertex orderings}.

Beyond describing connected puzzle assemblies, successive vertex orderings arise in many different contexts whenever a connected structure is built one component at a time. For example, in molecular self-assembly, they describe growth histories in which each newly attached molecule binds to an existing complex. In epidemic or information spreading on networks, they correspond to transmission histories in which each newly infected or informed individual is reached through an existing contact. 

Counting successive vertex orderings is a nontrivial combinatorial problem. Until recently, exact formulas for successive vertex orderings were known only for certain special families of graphs such as complete and complete bipartite graphs~\cite{fang2023}, but not for grid graphs arising from jigsaw puzzles. This limitation was overcome by the recent derivation of an exact counting formula for arbitrary finite connected graphs~\cite{agrawal2026successive}. In this article, we show how this framework can be used to enumerate and compare a variety of jigsaw-puzzle assembly strategies. As we shall see, the connectivity-preserving strategy that many people find most natural is actually quite rare: for larger puzzles, the overwhelming majority of assembly sequences pass through disconnected intermediate states before reaching the completed puzzle.

\section{A toy example:  counting assembly sequences of a \(3\times2\) puzzle}
\label{sec:direct_counting}

Before introducing the graph-theoretic framework, let us consider the six-piece jigsaw puzzle shown in Figure~\ref{fig:puzzle_to_graph}a, for which all $6!=720$ assembly sequences can be enumerated directly. 
We will consider three assembly strategies: connected growth from a single starting piece, connected growth from several disconnected starting pieces, and growth in which new disconnected components may be introduced at any stage during assembly. We use the term \emph{seed} exclusively for the initial piece, or group of disconnected pieces, placed at the very beginning of the puzzle. Any disconnected islands created later in the assembly are called new components.

\subsection{Single-seed connected assembly}

The most restrictive strategy is to insist that the partially assembled puzzle remains connected throughout the assembly process. Starting from a single piece, every newly placed piece must share an edge with at least one piece already placed. We call this strategy \emph{single-seed connected assembly} because it starts from a single seed piece and grows while remaining connected.

\begin{figure}
    \centering
    \includegraphics[width=0.8\linewidth]{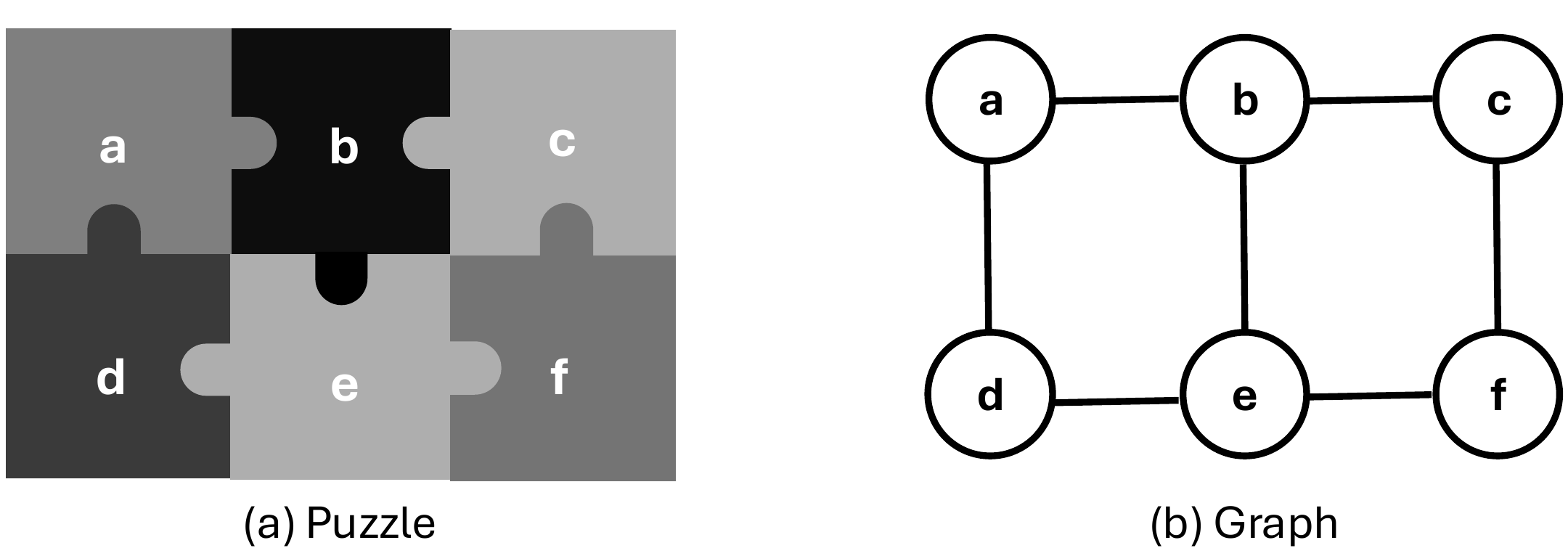}
    \caption{Jigsaw puzzle and its associated graph. (a) \(3\times2\) jigsaw puzzle with pieces labelled \(a\) through \(f\). (b) The corresponding graph $G(V,E)$ obtained by representing each puzzle piece as a vertex $ \in V$ and joining two vertices whenever the corresponding pieces share an edge $ \in E$. Assembly sequences of the puzzle can therefore be studied as orderings of the vertices of the associated graph.}
    \label{fig:puzzle_to_graph}
\end{figure}

For example, the ordering $(a \prec b \prec c \prec f \prec d \prec e)$ in Figure~\ref{fig:puzzle_to_graph} is a connected assembly sequence, where the symbol $\prec$ indicates that one piece is placed before the next. After each placement, the pieces already placed form a single connected cluster.
By contrast, \((a \prec d \prec c \prec f \prec b \prec e)\) is not connected, since the third piece placed, $c$, is disconnected from the preceding assembly.

By symmetry, the six pieces of the puzzle fall into two types: four corners and two middle pieces. If the first piece is a corner, there are two possible second pieces, and a short continuation count gives $24$ connected sequences starting from that corner. If the first piece is a middle piece, there are three possible second pieces, and the corresponding continuation count gives $56$ connected sequences starting from that piece. Hence the total number of connected assembly sequences is \(
4\cdot 24 + 2\cdot 56 = 208.\)

\subsection{Multi-seed connected assembly}

Many puzzle solvers do not begin from a single seed piece. Instead, they identify several promising pieces, such as the corners, place them first and then grow these disconnected regions before eventually merging them into the completed puzzle. In this strategy, after the initial seed pieces have been placed (in any order), every subsequent piece must attach to one of the existing partial assemblies; no new disconnected components are allowed to form. We call this strategy \emph{multi-seed connected assembly}. 

For example, the assembly sequence $(e \prec a \prec c \prec d \prec b \prec f)$ begins from three disconnected seed pieces.
Every subsequent piece attaches to one of the existing partial assemblies until the disconnected regions eventually merge into the completed puzzle.

To illustrate this assembly strategy, let us count the assembly sequences that start from three disconnected seed pieces. Since a $3\times2$ puzzle admits at most three mutually disconnected pieces, any newly added piece will inevitably attach to the existing assembly. There are only two possible choices of three pairwise disconnected seeds: $\{a,c,e\}$ and $\{b,d,f\}$. By symmetry, it suffices to consider the first of these.

Suppose the assembly begins from the seed set \(\{a,c,e\}\). The three seeds may be placed in any order, giving $3!$ possibilities. Once these seeds have been placed, the remaining pieces are \(b,d,\) and \(f\). Each of these pieces is already adjacent to at least one seed (denoted by the symbol $\sim$):
\(
b\sim \{a,c,e\},\,
d\sim \{a,e\},\,
f\sim \{c,e\}.
\)
Therefore every ordering of \(b,d,\) and \(f\) is valid, contributing a further $3!$ possibilities and  the seed set \(\{a,c,e\}\) contributes
\(
3!\times 3!=36
\)
assembly sequences. By symmetry, the seed set \(\{b,d,f\}\) contributes the same number. Thus,  the total number of assembly sequences that begin from three disconnected seeds and thereafter grow without creating any additional components is
\(
2\times 36 = 72.  
\)
Similar arguments can be used for two disconnected seed pieces. 

\subsection{Multi-component assembly}

The previous strategy allowed several disconnected seed pieces at the start, but required all subsequent growth to remain connected to one of the existing partial assemblies. We now consider a more general strategy in which new disconnected components may be introduced at any stage during assembly, including at the beginning, which subsumes multi-seed connected assembly.

In this strategy, the assembly begins from a single seed piece. At any subsequent stage of the assembly process, a newly placed piece may be disconnected from all previously placed pieces, thereby creating a new component. This island creation may happen immediately after the first piece is placed or much later during the assembly. Once a disconnected component is created, it continues to grow connectedly, with each newly added piece attaching to an existing partial assembly. Eventually, all disconnected components merge to produce the completed puzzle.
We call this strategy \emph{multi-component assembly}, since multiple disconnected components may be created at arbitrary stages of the assembly process.

For the puzzle pieces in Figure~\ref{fig:puzzle_to_graph}a, the ordering \(d \prec b \prec a \prec f \prec e \prec c\) begins with the seed piece \(d\). The second piece, \(b\), is not adjacent to \(d\), creating a new component. The third piece, \(a\), attaches to the existing assembly containing \(d\) and \(b\), so no new component is created. The fourth piece, \(f\), is disconnected and therefore creates another component. The remaining pieces, \(e\) and \(c\), are then added so as to connect and grow the existing components, eventually merging them into the completed puzzle.  

The total number of assembly sequences with one additional disconnected component is $424$ and with two additional disconnected components is $88$ (details in the Supplementary Material). Since there can be no more than three disconnected components, and this strategy subsumes the multi-seed connected assembly, the total number of assembly sequences is given by $208 + 424 + 88=720$.


\section{Representing jigsaw puzzles  as graphs and assembly sequences as vertex orderings }
\label{sec:puzzle_to_graphs}

The direct enumeration strategy used above for a simple $3 \times 2$ puzzle becomes infeasible for even slightly larger puzzles because the number of assembly paths grows as $n!$.
To make progress, we represent a jigsaw puzzle with $n$ pieces as a graph $G=(V,E)$ by mapping each puzzle piece to a vertex $v_i \in V$ and connecting two vertices by an edge whenever the corresponding pieces share a border. This defines a vertex set $V$ of cardinality $|V|=n$, and an edge set $E$. Figure~\ref{fig:puzzle_to_graph} illustrates this correspondence for the $3\times2$ puzzle.

Under this representation, every assembly sequence becomes a linear ordering of the vertices of the graph. The first vertex in the ordering represents the first piece placed, the second vertex represents the second piece placed, and so on. 
Let \(v_1 \prec v_2 \prec \cdots \prec v_n
\) be the vertex ordering corresponding to an assembly sequence. After the first \(k\) pieces have been placed, the partial assembly is represented by the subgraph induced by the vertices
\(
v_1,\ldots,v_k \in V
\).
The assembly strategies discussed above can therefore be reformulated as conditions on these growing subgraphs. For example, a \emph{single-seed connected assembly} sequence is one in which the vertices \(v_1,\ldots,v_k\) induce a connected subgraph for every \(1 \le k \le n\). In graph-theoretic language, such orderings are known as \emph{successive vertex orderings}. 

The second assembly strategy, \emph{multi-seed connected assembly}, begins from \(s\) seed pieces, corresponding to the vertices \(v_1,\ldots,v_s\), where \(s<n\). These vertices are not adjacent to each other and thus form an independent set. For every \(k>s\), however, the vertex \(v_k\) must have a neighbour among \(v_1,\ldots,v_{k-1}\). Consequently, the initial partial assemblies can grow and eventually merge together, but no additional disconnected components can be created after the seed set.

The third assembly strategy, \emph{multi-component assembly}, allows additional disconnected components to be introduced after placing the seed piece. Specifically, for some \(1<k<n\), the vertex \(v_k\) may be placed before all of its neighbours. Each such vertex initiates a new disconnected component. More generally, an assembly sequence is said to introduce \(k\) additional disconnected components if exactly \(k\) vertices, other than the initial seed vertex, are placed before all of their neighbours.

\section{Counting jigsaw puzzle assembly sequences using graph theory}
\label{sec:counting}

Having formulated the three assembly strategies in graph-theoretic language, we now show how the mathematical formalism for counting vertex orderings developed in \cite{fang2023,agrawal2026successive} can be applied to jigsaw puzzles. We  state the relevant counting formula corresponding to each strategy, briefly sketching out the intuition behind the fuller mathematical treatment,  and illustrate their use on the \(3\times2\) puzzle shown in Figure~\ref{fig:puzzle_to_graph}. 


\subsection{Single-seed connected assembly}

A single-seed connected assembly fails precisely when a newly placed piece does not attach to the existing assembly, forming a new component. In our graph representation, this occurs when a non-first vertex appears in the ordering before all of its neighbours. In \cite{fang2023,agrawal2026successive}, such a vertex is called a \emph{bad} vertex. Therefore, a vertex ordering corresponds to a single-seed connected assembly sequence if and only if it contains no bad vertices.

Let $\sigma(G)$ denote the total number of successive vertex orderings. Since there are $n!$ possible vertex orderings in total, the ratio $\sigma(G)/n!$ can be interpreted as the probability that a randomly and uniformly chosen linear vertex ordering  contains no bad vertices. We use the language of probabilities because it makes the analysis easier.     For any subset of vertices $I \subseteq V$, let $B_I$ denote the event that every vertex in $I$ is bad simultaneously. The inclusion-exclusion principle of set theory can be used to write  the probability of a successive vertex ordering of a graph $G$, 
\begin{equation}
\label{eq:inc-exc}
   \frac{\sigma(G)}{n!} = \sum_{I \subseteq V} (-1)^{|I|} \mathbb{P}(B_I), 
\end{equation} 
as an alternating sum~(\ref{eq:inc-exc}) over all possible subsets $I \subseteq V$  of the probability $\mathbb{P}(B_I)$ that all vertices in $I$ are bad, in a random ordering.

Most subsets contribute nothing to this sum. If $I$ contains two adjacent vertices, it is logically impossible for both to be bad, as each would have to appear in the ordering before the other. Therefore, $\mathbb{P}(B_I) = 0$ whenever $I$ contains adjacent vertices. Consequently, only  independent sets, i.e.\ sets in which no two vertices are adjacent,  contribute non-zero terms to the sum.

To calculate $\mathbb{P}(B_I)$, we decompose the ordering process into two logical steps:

\begin{itemize}
    \item {The first vertex:} To determine which vertex may appear first in the ordering while ensuring that every vertex of $I$ remains bad, let $N[I]$ denote the closed neighbourhood of $I$, consisting of the vertices of $I$ together with all vertices adjacent to at least one vertex of $I$. If the first vertex in the ordering belonged to $N[I]$, then at least one vertex of $I$ would immediately acquire an earlier neighbour, failing the condition. Hence, the first vertex must lie strictly outside $N[I]$. The number of valid choices for the first vertex is therefore given by the quantity:
    \begin{equation}
    \label{eq:aI}
    a(I)=|V\setminus N[I]|=n-|N[I]|.
    \end{equation}
    Since every vertex is equally likely to appear first in a uniformly random ordering, the probability of making a valid first choice is exactly $a(I)/n$.
    
    \item {The remaining vertices:} After placing the first vertex, we must ensure that all vertices in $I$ appear in the ordering before any of their neighbours. We define $b(I)$ as the conditional probability that this ordering succeeds. This probability follows a recursive structure:
    \begin{equation}
    \label{eq:bI}
    b(\emptyset)=1, \qquad b(I)=\frac{1}{|N[I]|} \sum_{v\in I} b(I\setminus\{v\}), \qquad I\neq\emptyset,
    \end{equation}
    where $|N[I]|$ is the size of the closed neighbourhood of $I$. This recursion has a precise logical interpretation. Among the vertices of $I$, one of them must be the first to appear in the ordering. Suppose this vertex is $v$. As the earliest vertex in $I$, $v$ precedes all other vertices in the set. Since every vertex in $I$ is bad and must therefore precede its own neighbours, $v$ logically precedes all of them as well, making $v$ the earliest vertex in the entire closed neighbourhood $N[I]$. The probability of this event is $1/|N[I]|$. Once $v$ has appeared, the remaining vertices of $I\setminus\{v\}$ must satisfy exactly the same conditions. Their contribution is therefore $b(I\setminus\{v\})$. Summing over every possible choice for this first bad vertex produces the recursion above.
    
\end{itemize}

Combining the probability of choosing a valid first vertex, $a(I)/n$, with the factor $b(I)$ gives exactly $\mathbb{P}(B_I) = a(I)b(I)/n$.  Substituting this into our inclusion-exclusion sum and multiplying by $n!$ yields the exact formula for the number of successive vertex orderings $\sigma(G)$~\cite[Theorem~1.2]{agrawal2026successive}:
\begin{equation}
\label{eq:sigma_G}
\sigma(G) = n!\sum_{\substack{I\subseteq V\\ I\text{ independent}}} (-1)^{|I|} \frac{a(I)}{n}\,b(I).
\end{equation}

 
\subsubsection{Example} To illustrate how this calculation works, we apply Eq.~\eqref{eq:sigma_G} to the grid graph in Figure~\ref{fig:puzzle_to_graph}b where \(n=6\). The independent sets fall into a small number of symmetry classes, and the corresponding values of \(a(I)\) and \(b(I)\) are shown in Table~\ref{Table:independent_sets_3x2}.
\begin{table}[b]
\centering
\caption{Independent sets of the \(3\times2\) puzzle's graph grouped into symmetry classes, together with the corresponding values of \(a(I)\) and \(b(I)\).}
\label{Table:independent_sets_3x2}
\renewcommand{\arraystretch}{1.15}
\begin{tabular}{lccc}
\hline
Representative independent set(s) & multiplicity & \(a(I)\) & \(b(I)\) \\
\hline
\(\emptyset\) & \(1\) & \(6\) & \(1\) \\[2pt]
\(\{a\},\{c\},\{d\},\{f\}\) & \(4\) & \(3\) & \(1/3\) \\
\(\{b\},\{e\}\) & \(2\) & \(2\) & \(1/4\) \\[2pt]
\(\{a,e\},\{b,d\},\{b,f\},\{c,e\}\) & \(4\) & \(1\) & \({7}/{60}\) \\
\(\{a,c\},\{d,f\}\) & \(2\) & \(1\) & \({2}/{15}\) \\
\(\{a,f\},\{c,d\}\) & \(2\) & \(0\) & \(1/9\) \\[2pt]
\(\{a,c,e\},\{b,d,f\}\) & \(2\) & \(0\) & \({11}/{180}\) \\
\hline
\end{tabular}
\end{table}
These values are obtained directly from the definition of \(a(I)\) and the recursion for \(b(I)\). For example,
\[
a(\{a\})=|V\setminus N[\{a\}]|=3,
\qquad
b(\{a\})=\frac{1}{|N[\{a\}]|}=\frac13,
\]
and
\[
b(\{a,e\})
=
\frac{1}{|N[\{a,e\}]|}\bigl(b(\{a\})+b(\{e\})\bigr)
=
\frac{1}{5}\left(\frac13+\frac14\right)
=
\frac{7}{60}.
\]
Similarly,
\[
b(\{a,c\})
=
\frac{1}{5}\left(\frac13+\frac13\right)
=
\frac{2}{15},
\qquad
b(\{a,c,e\})
=
\frac{1}{6}\left(\frac{7}{60}+\frac{2}{15}+\frac{7}{60}\right)
=
\frac{11}{180}.
\]

Substituting these values into the formula for \(\sigma(G)\) gives
\[
\frac{\sigma(G)}{6!}
=
1
- \left(4\cdot \frac{3}{6}\cdot \frac13 + 2\cdot \frac{2}{6}\cdot \frac14\right)
+ \left(4\cdot \frac{1}{6}\cdot \frac{7}{60} + 2\cdot \frac{1}{6}\cdot \frac{2}{15}\right),
\]
since the terms with \(a(I)=0\) do not contribute. Hence
\[
\frac{\sigma(G)}{6!}
=
1-\frac56+\frac{11}{90}
=
\frac{13}{45} \implies \sigma(G)=6!\cdot \frac{13}{45}=208.
\]
This result agrees with the value obtained earlier by direct counting.  

\subsection{Multi-seed connected assembly}

Let $\tau(G,S)$ count the connected assembly sequences of $G$ that start from a prescribed independent seed set $S$. After the seeds are placed, each new vertex must attach to the partial assembly, and this can only fail at a vertex with no neighbour in $S$. Setting $H = G - S$, these are exactly the vertices of $W = V \setminus N[S]$, the complement of the closed neighbourhood of the seeds. The argument then runs exactly as in the single-seed case: inclusion--exclusion over the independent sets of $H$ lying inside $W$, starting from the probability that every vertex of an independent set $I \subseteq W$ is bad.

Since the first vertex is already selected from $S$, we no longer need the factor $a(I)$ to calculate $\mathbb P(B_I)$.
Consequently, the probability that every vertex in $I$ is bad is given entirely by the correction factor $b_H(I)$. This factor follows the exact same recursion as before, except all neighbourhoods are now evaluated within the reduced subgraph $H = G - S$:
\begin{equation}
\label{eq:bH}
b_H(\emptyset)=1, \qquad
b_H(I) = \frac{1}{|N_H[I]|} \sum_{v\in I} b_H(I\setminus\{v\}), \qquad I\neq\emptyset,
\end{equation}
where $N_H[I]$ is the closed neighbourhood of $I$ within $H$.

To convert this probability into an absolute count of orderings, we multiply it by the total number of unconstrained arrangements allowed by this strategy. There are $|S|!$ ways to arbitrarily order the initial seeds, and $(n-|S|)!$ ways to order the remaining vertices. Multiplying the probability $b_H(I)$ by these two factors gives the exact number of orderings in which $I$ is bad. Substituting this into the alternating inclusion-exclusion sum in Eq.~\eqref{eq:inc-exc} yields the final formula \cite[Theorem B.2]{agrawal2026successive}:
\begin{equation}
\label{eq:tauGS}
\tau(G,S) = |S|!(n-|S|)! \sum_{\substack{I\subseteq W\\ I\text{ independent}}} (-1)^{|I|} \,b_H(I).
\end{equation}

\subsubsection{Example}
As an illustrative example, we compute the number of assembly sequences that begin from three disconnected seeds and thereafter grow without creating any additional disconnected components. As shown previously, there are only two possible choices of three pairwise disconnected seeds:
\(
\{a,c,e\}\text{ and } \{b,d,f\}.
\)
By symmetry, it suffices to consider \(S=\{a,c,e\}\).

Every remaining vertex \(b,d,f\) is adjacent to at least one vertex of \(S\). Consequently,
\[
W=V\setminus N[S]=\varnothing.
\]
The only independent set contained in \(W\) is therefore the empty set, for which \(b_H(\emptyset)=1\). Eq.~\eqref{eq:tauGS} then gives
\[
\tau(G,\{a,c,e\})
=
3!\,3!
=
36.
\]

By symmetry,
\(
\tau(G,\{b,d,f\})=36.
\)
Hence the total number of assembly sequences that begin from three disconnected seed pieces and thereafter grow without introducing any additional disconnected components is $36+36=72$, which agrees  with the value obtained earlier by direct combinatorial counting.

One example of where this formalism can be used is for the popular strategy of identifying the four corner pieces and placing them first.
For an \(M\times N\) puzzle (where \(M,N\ge3\)), these four vertices form the independent seed set
 \(
 S=\{\text{four corner vertices}\}
 \). Once the corner pieces have been fixed, every subsequent piece must attach to one of the existing partial assemblies. The number of such assembly sequences is therefore obtained directly from Eq.~\eqref{eq:tauGS}.



\subsection{Multi-component assembly}

This strategy permits any non-first vertex to appear in the ordering before all of its neighbours. Each such bad vertex creates one new disconnected component. Rather than counting only successive vertex orderings, which are precisely the orderings with no bad vertices, we now wish to enumerate orderings according to the number of bad vertices they contain.

Let \(A_k\) denote the number of vertex orderings with exactly \(k\) bad vertices. Equivalently, \(A_k\) is the number of assembly sequences in which exactly \(k\) additional disconnected components are introduced after placing the initial seed piece. 

The inclusion-exclusion formula for connected assembly in Eq.~\eqref{eq:sigma_G} already encodes information about bad vertices. However, the contributions from different independent sets are immediately collapsed by the alternating factor \((-1)^{|I|}\), leaving only the final count \(A_0=\sigma(G)\). To preserve the information carried by the individual terms, \cite{agrawal2026successive} replaced the alternating factor with a variable \(x^{|I|}\).

This substitution defines a polynomial over the independent sets of the graph. It is structurally similar to the classic independence polynomial in graph theory, but each term is now weighted by the  factor \(a(I)b(I)/n\). This weighted sum is known as the \emph{successive ordering polynomial}:
\begin{equation}
\label{eq:polynomial}
P_G(x) = \sum_{\substack{I\subseteq V\\I\text{ independent}}} \frac{a(I)}{n}\,b(I)\,x^{|I|}.
\end{equation}
The variable \(x\) in Eq.~\eqref{eq:polynomial} is best understood as a marker for chosen bad vertices. The exponent of \(x\) records the size of the independent set \(I\). But the  factor \(a(I)b(I)/n\) attached to \(x^{|I|}\) has a more specific meaning: it is the probability that every vertex in \(I\) is bad. Thus a term of the polynomial records orderings in which a chosen set \(I\) of vertices is bad, and the power of \(x\) records how many such bad vertices have been chosen.

The important subtlety is that the event ``all vertices in \(I\) are bad'' does not say that these are the only bad vertices. If a particular ordering has several bad vertices, then it is counted once for every subset of its bad vertices. For example, suppose an ordering has exactly three bad vertices, say \(\{u,v,w\}\). There is one way to choose none of them, three ways to choose one of them, three ways to choose two of them, and one way to choose all three. Since choosing \(j\) bad vertices contributes \(x^j\), this ordering contributes
\(
1+3x+3x^2+x^3=(1+x)^3.
\)
More generally, if an ordering has \(r\) bad vertices, then choosing \(j\) of them can be done in \(\binom{r}{j}\) ways, so its total contribution is
\(
\sum_{j=0}^{r}\binom{r}{j}x^j=(1+x)^r.
\)

Now consider all orderings with exactly \(r\) bad vertices. Each such ordering contributes \((1+x)^r\), and by definition there are \(A_r\) such orderings. Therefore, the total contribution from all orderings with exactly \(r\) bad vertices is
\(
A_r(1+x)^r.
\)
Adding these contributions over all possible values of \(r\) gives
\[
n!P_G(x)
=
A_0
+
A_1(1+x)
+
A_2(1+x)^2
+
A_3(1+x)^3
+\cdots .
\]


We now need a way to extract the coefficient of \((1+x)^k\). This is exactly what differentiating at \(x=-1\) does. At \(x=-1\), the factor \(1+x\) becomes zero. Evaluating \(P_G(x)\) at \(x=-1\) therefore keeps only the constant term in the expansion above, namely the term with no bad vertices:
\[
n!P_G(-1)=A_0=\sigma(G).
\]
To extract the next term, we differentiate first. Differentiating removes one factor of \(1+x\). Thus the term \(A_1(1+x)\) becomes \(A_1\), while every term \(A_r(1+x)^r\) with \(r\ge 2\) still contains at least one remaining factor of \(1+x\), and therefore vanishes when \(x=-1\). Hence
\[
n!P_G'(-1)=A_1.
\]
Similarly,
\[
\frac{n!}{2!
}P_G''(-1)=A_2.
\]
The same reasoning gives, for every \(k\) \cite[Theorem~5.3]{agrawal2026successive},
\begin{equation}
\label{eq:Ak}
    A_k=\frac{n!}{k!}P_G^{(k)}(-1).
\end{equation}

This single formula unifies all three assembly strategies. When \(k=0\), we recover \(A_0=\sigma(G)\), corresponding to single-seed connected assembly. If several components are introduced at the very beginning and no new components are allowed to appear later, the assembly follows the multi-seed connected strategy. The multi-component strategy is the most general one, since it allows these additional components to appear at arbitrary stages during assembly.

\subsubsection{Example}

For the \(3\times2\) puzzle, Eq.~\eqref{eq:polynomial} yields
\[
P_G(x)
=
1+\frac{5}{6}x+\frac{11}{90}x^2.
\]

Solving for $A_k$ using Eq.~\eqref{eq:Ak} gives
\[
A_0
=
6!\,P_G(-1)
=
208; \quad
A_1
=
6!\,P_G'(-1)
=
424; \quad
A_2
=
\frac{6!}{2}\,P_G''(-1)
=
88.
\]
Thus there are \(208\) assembly sequences that remain connected throughout, \(424\) assembly sequences in which one additional disconnected component is introduced, and \(88\) assembly sequences in which two additional disconnected components are introduced.
The values agree with the direct counts obtained earlier. This framework recovers all the counts of different assembly strategies without requiring a case-by-case analysis.

\begin{table}[b]
\caption{
Number of multi-component assembly sequences for \(3\times3\), \(4\times4\), \(5\times5\), and \(6\times6\) jigsaw puzzles. The quantity \(A_k\) counts assembly sequences in which exactly \(k\) disconnected components are created after the initial seed piece is placed.
}
\centering
\small
\begin{tabular}{ccccc}
\hline
\textbf{Number of disconnected} 
&
\textbf{\(A_k\) for \(3\times3\)}
&
\textbf{\(A_k\) for \(4\times4\)}
&
\textbf{\(A_k\) for \(5\times5\)}
&
\textbf{\(A_k\) for \(6\times6\)} \\
\textbf{components \(k\)}
& & & & \\
\hline

0 
& \(2.96 \times 10^{4}\)
& \(3.57 \times 10^{10}\)
& \(8.84 \times 10^{19}\)
& \(1.08 \times 10^{33}\) \\
1 
& \(1.44 \times 10^{5}\)
& \(7.79 \times 10^{11}\)
& \(7.06 \times 10^{21}\)
& \(3.04 \times 10^{35}\) \\
2 
& \(1.61 \times 10^{5}\)
& \(4.65 \times 10^{12}\)
& \(1.56 \times 10^{23}\)
& \(2.21 \times 10^{37}\) \\
3 
& \(2.40 \times 10^{4}\)
& \(8.31 \times 10^{12}\)
& \(1.15 \times 10^{24}\)
& \(5.10 \times 10^{38}\) \\
4 
& \(4.80 \times 10^{3}\)
& \(5.39 \times 10^{12}\)
& \(3.46 \times 10^{24}\)
& \(4.88 \times 10^{39}\) \\
5
& 
& \(1.58 \times 10^{12}\)
& \(4.97 \times 10^{24}\)
& \(2.33 \times 10^{40}\) \\
6 
& 
& \(1.67 \times 10^{11}\)
& \(3.74 \times 10^{24}\)
& \(6.14 \times 10^{40}\) \\
7 
& 
& \(1.44 \times 10^{10}\)
& \(1.58 \times 10^{24}\)
& \(9.60 \times 10^{40}\) \\
8 
& 
& 
& \(3.85 \times 10^{23}\)
& \(9.33 \times 10^{40}\) \\
9 
& 
& 
& \(6.04 \times 10^{22}\)
& \(5.87 \times 10^{40}\) \\
10 
& 
& 
& \(6.19 \times 10^{21}\)
& \(2.48 \times 10^{40}\) \\
11 
& 
& 
& \(6.22 \times 10^{20}\)
& \(7.25 \times 10^{39}\) \\
12 
& 
& 
& \(4.78 \times 10^{19}\)
& \(1.50 \times 10^{39}\) \\
13 
& 
& 
& 
& \(2.24 \times 10^{38}\) \\
14 
& 
& 
& 
& \(2.47 \times 10^{37}\) \\
15 
& 
& 
& 
& \(2.08 \times 10^{36}\) \\
16 
& 
& 
& 
& \(1.51 \times 10^{35}\) \\
17 
& 
& 
& 
& \(6.97 \times 10^{33}\) \\
\hline
\end{tabular}
\label{tab:Ak_comparison}
\end{table}

\begin{figure}
    \centering
    \includegraphics[width=0.8\linewidth]{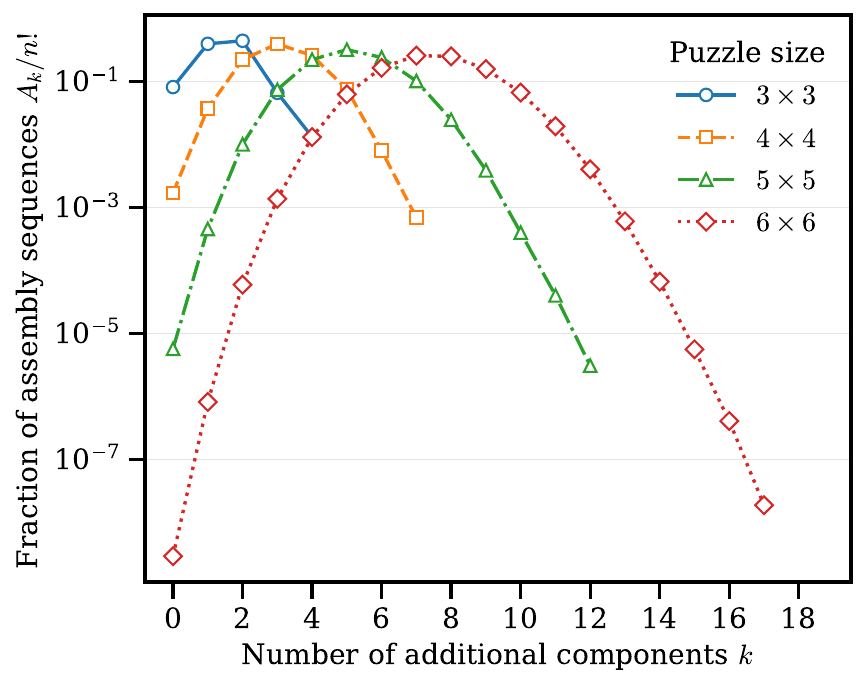}
    \caption{Distribution of assembly sequences by the number of additional disconnected components introduced during assembly. Each curve shows the fraction $A_k/n!$ of assembly sequences for an $n = N\times N$ jigsaw puzzle in which exactly $k$ additional disconnected components are created after the initial seed piece. For a fixed puzzle size, assembly sequences are most likely to create an intermediate number of disconnected components, whereas sequences that remain nearly connected or create very many disconnected components occur much less frequently.}
    \label{fig:Ak_dist}
\end{figure}

\begin{table}
\caption{Comparison of two assembly strategies for a \(5\times5\) jigsaw puzzle that each involve \(k\) disconnected components. The quantity \(\sum_{\substack{|S|=k+1,\, S\text{ ind.}}}\tau(G,S)\) counts assembly sequences that begin from \(k+1\) disconnected seeds and thereafter remain connectivity-preserving. In contrast, \(A_{k}\) counts assembly sequences that begin from a single seed and later introduce \(k\) additional disconnected components at arbitrary stages during assembly.}
\begin{center}
\small
\begin{tabular}{ccc}
\hline
\begin{tabular}{c}
\textbf{Number of
disconnected}\\ \textbf{components}\\
\(k\)
\end{tabular}
&
\begin{tabular}{c}
\textbf{Introduced as seeds}\\
\(\displaystyle \sum_{\substack{|S|=k+1,\\S\text{ independent}}}\tau(G,S)\)
\end{tabular}
&
\begin{tabular}{c}
\textbf{Introduced at arbitrary} \\\textbf{stages during assembly}\\
\(A_{k}\)
\end{tabular}
\\
\hline
\(0\)  & \(8.84 \times 10^{19}\) & \(8.84 \times 10^{19}\) \\
\(1\)  & \(2.91 \times 10^{21}\) & \(7.06 \times 10^{21}\) \\
\(2\)  & \(2.90\times 10^{22}\) & \(1.56 \times 10^{23}\) \\
\(3\)  & \(1.07 \times 10^{23}\) & \(1.15 \times 10^{24}\) \\
\(4\)  & \(1.89 \times 10^{23}\) & \(3.46 \times 10^{24}\) \\
\(5\)  & \(1.86 \times 10^{23}\) & \(4.97 \times 10^{24}\) \\
\(6\)  & \(1.10 \times 10^{23}\) & \(3.74 \times 10^{24}\) \\
\(7\)  & \(4.18 \times 10^{22}\) & \(1.58 \times 10^{24}\) \\
\(8\)  & \(1.04 \times 10^{22}\) & \(3.85 \times 10^{23}\) \\
\(9\) & \(1.81 \times 10^{21}\) & \(6.04 \times 10^{22}\) \\
\(10\) & \(2.39 \times 10^{20}\) & \(6.19 \times 10^{21}\) \\
\(11\) & \(3.18 \times 10^{19}\) & \(6.22 \times 10^{20}\) \\
\(12\) & \(2.98 \times 10^{18}\) & \(4.78 \times 10^{19}\) \\
\hline
\end{tabular}
\end{center}
\label{Table:Ak_tau}
\end{table}

\section{Assembly statistics for $N \times N$ puzzles}
\label{sec:larger_puzzles}

Having established the exact algebraic formulas for our assembly strategies, we can now apply this graph-theoretic framework to compute assembly statistics for puzzles far beyond the reach of direct enumeration. As an illustration, consider the $5\times5$ puzzle shown in Figure~\ref{fig:mona_lisa}a.
Enumerating the independent sets of the corresponding grid graph and evaluating the associated quantities $a(I)$ and $b(I)$ gives the successive ordering polynomial:
\[
\begin{aligned}
P_G(x) ={}& 1
+ \frac{77}{15}x
+ \frac{239}{21}x^2
+ \frac{23222851}{1621620}x^3
+ \frac{6232148353}{550368000}x^4 \\
&+ \frac{2981890237939}{505237824000}x^5 
+ \frac{3254514069009079}{1565286218496000}x^6 \\
&+ \frac{737376332468414501}{1463542614293760000}x^7 
+ \frac{201286386772326871}{2341668182870016000}x^8\\
&+ \frac{1432883497675129}{133049328572160000}x^9 
+ \frac{537977544649823}{515391083311104000}x^{10}\\
&+ \frac{3322019035145479}{43086694564808294400}x^{11} 
+ \frac{3322019035145479}{1077167364120207360000}x^{12}.
\end{aligned}
\]

Applying Eq.~\eqref{eq:Ak} to this polynomial extracts $A_k$ for the $5\times5$ puzzle. Similar polynomials can be calculated for other puzzles, and the complete distribution of $A_k$ across $3\times3$, $4\times4$, $5\times5$, and $6\times6$ puzzles is summarized in Table~\ref{tab:Ak_comparison}. (The code used to generate these results is available at \cite{svocode2026}.) 
Since the largest possible number of disconnected components is given by the independence number, defined as  the size of the largest independent set  of the corresponding grid graph, the table terminates at different values of \(k\) for different puzzle sizes. Since the first piece placed acts as the initial seed, the number of additional components $k$ cannot exceed the independence number minus one. For example, the \(5\times5\) grid has independence number \(13\), so \(k\le12\).

In Figure~\ref{fig:Ak_dist} we show  the normalized distributions \(A_k/n!\) from Table~\ref{tab:Ak_comparison}.
A number of clear patterns emerge from this table and figure.   First, growing a single connected cluster, which may be the most intuitive assembly structure, only accounts for a tiny fraction of all assembly sequences. Moreover, this fraction rapidly decreases with increasing puzzle size. 

Second, the distributions are peaked at intermediate values of \(k\). The most numerous assembly sequences introduce a moderate number of additional disconnected components, whereas sequences that introduce either very few or very many are comparatively rare. This behaviour reflects a combinatorial trade-off. Requiring the assembly to remain nearly connected severely restricts the allowable orderings, while creating many disconnected components requires many pieces to be placed before all of their neighbours, which is itself increasingly restrictive. The competition between these opposing constraints produces the broad peaks observed in Figure~\ref{fig:Ak_dist}.

We can also compare two ways of reaching the same number of disconnected components: introducing them all at the start, as in multi-seed connected assembly, or letting them appear at arbitrary stages, as counted by $A_k$. 
Since every multi-seed connected assembly is also a multi-component assembly, the corresponding counts satisfy
\[
\sum_{\substack{|S|=k+1\\S\text{ independent}}}
\tau(G,S)
\le
A_k,
\]
with equality only when no additional disconnected components are introduced after the initial seed pieces.

Table~\ref{Table:Ak_tau} illustrates this comparison for the \(5\times5\) puzzle. The second column gives the total number of assembly sequences that begin from \(k+1\) disconnected seed pieces and thereafter remain connectivity-preserving, summed over all independent seed sets. The third column gives the number of assembly sequences that begin from a single seed and introduce \(k\) additional disconnected components at arbitrary stages during assembly.

The two counts coincide when \(k=0\), corresponding to connected assembly. For every \(k>0\), however, the multi-component counts are substantially larger than the corresponding multi-seed connected counts. Thus, allowing disconnected components to appear during the assembly process dramatically increases the number of possible assembly sequences. Most assembly sequences with \(k+1\) disconnected components do not create all of those components at the outset; instead, new disconnected regions typically arise progressively as the puzzle is assembled.

This formalism allows us to compare different starting strategies for connected growth.   Consider the $5\times 5$ puzzle. If you put down one fixed corner piece, the total number of connected sequences is  $6.95 \times 10^{16}$.  The popular strategy of placing four corners at the beginning gives  $9.03 \times 10^{18}$ connected sequences, which is considerably larger because you can choose which of the four assemblies to add a new piece to.   If you start with the centre piece, you have even more connected assembly sequences available, $1.76 \times 10^{19}$.  It's easy to see why this gives more assembly sequences than a single corner piece.  The reason this gives more assemblies than the four corners is because in the latter strategy, there are $21$ pieces left to place, whereas for a single seed, there are $24$.   Many more strategies can be compared in this way. 


\section{Conclusion}

We began by asking how many ways a jigsaw puzzle can be assembled one piece at a time while remaining connected throughout the assembly process. We then extended this question to two broader assembly strategies: one in which assembly begins from a set of disconnected seed pieces and thereafter remains connectivity-preserving, and another in which new disconnected components may be introduced at arbitrary stages before eventually merging into the completed puzzle. 

By representing a jigsaw puzzle as a graph, each of these questions becomes a problem of counting vertex orderings, with recent results in graph theory \cite{agrawal2026successive} allowing for exact enumeration.  However, the step of enumerating all independent sets is slow, limiting this approach to relatively small puzzles.  

The resulting counts reveal that assembly sequences where the pieces remain connected throughout are extremely rare. For larger puzzles, the overwhelming majority of assembly sequences pass through disconnected intermediate stages, with the most common sequences creating an intermediate number of disconnected components rather than very few or very many.

Although motivated here by jigsaw puzzles, the quantities
$\sigma(G)$, $\tau(G,S)$, and $A_k$ apply to any process in which a connected structure is built one piece at a time.
Examples range from molecular self-assembly and crystal growth to assembly planning in robotics and manufacturing, where one may wish either to preserve connectivity throughout construction or to allow multiple disconnected regions that later merge.
Viewed in this way, a jigsaw puzzle provides a surprisingly rich model for studying incremental growth on graphs.

Many  questions remain open. For example: 
Is there a useful asymptotic approximation for
$\sigma(G)$  as  puzzle size grows? Similarly, does the simple form of the distribution of  $A_k$  in~Figure~\ref{fig:Ak_dist} suggest  an approximate method to allow extensions to larger puzzles?   Even more ambitiously, is there a way to exploit the symmetry of grid graphs to simplify the direct calculation of $\sigma(G)$, allowing solutions for larger puzzles? 
Which disconnected seed configurations maximize $\tau(G,S)$?
How do these quantities depend on the geometry of the underlying graph, and what changes for graphs that are far from grid-like? 

\bibliographystyle{vancouver}
\bibliography{references}

\clearpage
\setcounter{section}{0}
\renewcommand{\thesection}{A}
\renewcommand{\thefigure}{A\arabic{figure}}
\setcounter{figure}{0}
\renewcommand{\thetable}{A\arabic{table}}
\setcounter{table}{0}
\renewcommand{\theequation}{A\arabic{equation}}
\setcounter{equation}{0}

\section*{\textbf{Supplementary Material}}

In this supplementary material, we provide explicit enumeration details for the assembly-sequence counts reported in the main text. The goal is to make clear how assembly sequences are classified according to the number of additional disconnected components introduced during the assembly process. We focus on the
smallest nontrivial rectangular puzzle, the \(3\times 2\) grid, where all possible assembly sequences can be counted directly. 

We organize the count in the following order, starting with the case of two additional
disconnected components, then treating one additional disconnected component,
and finally the connected case with zero additional disconnected components.

\section{Counting assembly sequences with two disconnected components}

Here, we directly count the number of assembly sequences in which exactly two additional disconnected components are introduced after the seed vertex. By the symmetry of the underlying grid graph, the six starting pieces fall
into two types: \emph{corner} ($a,c,d,f$) and \emph{middle} ($b,e$).
We fix the first piece in each case and count sequences of the remaining five pieces that contain exactly two additional disconnected pieces. Recall that a piece is called 'bad' is it is not placed first and it appears before all of its neighbours.

\medskip\noindent
\textbf{Corner start} (first piece $= a$, neighbours $b,d$).

Once $a$ is placed, $b$ and $d$ each already have a placed neighbour
($a$ itself), so they can never be bad.  
The only pieces that can be bad are:
\begin{itemize}
\item $c$, if it is placed before both $b$ and $f$;
\item $e$, if it is placed before $b$, $d$, and $f$;
\item $f$, if it is placed before both $c$ and $e$.
\end{itemize}

Note that $e$ bad requires $e\prec f$, while $f$ bad requires $f\prec e$;
these conditions are mutually exclusive, so $e$ and $f$ cannot both be bad.
Similarly $c$ bad requires $c\prec f$, while $f$ bad requires $f\prec c$,
so $c$ and $f$ cannot both be bad.  Hence the only pair that can
simultaneously be bad is $\{c,e\}$.

We count orderings of $\{b,c,d,e,f\}$ in which both $c$ and $e$ are bad:
$e\prec\{b,d,f\}$ and $c\prec\{b,f\}$, which together require $b$ and
$f$ to follow \emph{both} $c$ and $e$.

Enumerate by the position of $e$ in the ordering (since $c$ and $e$ must both precede
$b$ and $f$, neither $b$ nor $f$ can appear in the first two positions after $a$) as in Table~\ref{tab:corner_count}:

\begin{table}[h]
\centering
\caption{Counting orderings in the corner-start case where both \(c\) and \(e\) are bad.}
\label{tab:corner_count}
\renewcommand{\arraystretch}{1.2}
\begin{tabular}{c p{6cm} p{3.5cm}}
\toprule
Position of $e$ &
Relative arrangements of $\{b,c,e,f\}$ &
Ways to insert $d$ (after $e$) \\
\midrule
2nd &
$e$ first, $c$ second, then $b,f$ in $2!$ orders: $2$ arrangements &
$d$ in positions $3$--$6$: $4$ ways \\
3rd &
$c$ first, $e$ second, then $b,f$ in $2!$ orders: $2$ arrangements &
$d$ in positions $4$--$6$: $3$ ways \\
\bottomrule
\end{tabular}
\end{table}

Total: $(2\times 4)+(2\times 3)=14$ orderings.

One corner piece contributes $14$ sequences; four corner pieces give
$4\times 14=56$.

\medskip\noindent
\textbf{Middle start} (first piece $= b$, neighbours $a,c,e$).

Once $b$ is placed, $a$, $c$ and $e$ each already have a placed neighbour
($b$ itself), so they are never bad. The only remaining pieces that
can be bad are:
\begin{itemize}
    \item $d$, if it is placed before both $a$ and $e$;
    \item $f$, if it is placed before both $c$ and $e$.
\end{itemize}

For both $d$ and $f$ to be bad:
$d\prec\{a,e\}$ and $f\prec\{c,e\}$, which together require $e$ to follow both $d$ and $f$.

Enumerate orderings of $\{a,c,d,e,f\}$ by the position of $e$ (which
must be at least fourth, since both $d$ and $f$ precede it) as in Table~\ref{tab:middle_start}:
\begin{table}[h]
\caption{Counting orderings in the middle-start case by the position of \(e\), assuming both \(d\) and \(f\) are bad.}
\label{tab:middle_start}
\begin{center}
\renewcommand{\arraystretch}{1.2}
\begin{tabular}{c p{8.5cm} c}
\hline
Position of $e$ & Remaining ordering constraints & Count \\
\hline
4th &
Positions \(2,3\) must contain \(\{d,f\}\) in any order (\(2!\) ways),
while positions \(5,6\) contain \(\{a,c\}\) in any order (\(2!\) ways).
&
\(2\times2=4\)
\\
5th &
Position \(6\) is either \(a\) or \(c\). Positions \(2,3,4\)
contain the remaining pieces, subject to either \(d\prec a\) or
\(f\prec c\), giving \(3!/2=3\) orderings in each case.
&
\(3+3=6\)
\\
6th &
Positions \(2\)--\(5\) contain \(\{a,c,d,f\}\) with the independent
constraints \(d\prec a\) and \(f\prec c\). Hence the number of
orderings is \(4!/4=6\).
&
\(6\)
\\
\hline
\multicolumn{2}{r}{Total} & \(16\) \\
\hline
\end{tabular}
\end{center}
\end{table}
One middle piece contributes $16$ sequences; two middle pieces give $2\times 16=32$. Combining both symmetry classes gives a total of \(56+32=88\) assembly sequences with exactly two additional disconnected components.

\section{Counting assembly sequences with one disconnected component}

We next count the assembly sequences in which exactly one additional
disconnected piece is introduced after the seed vertex. 

\medskip\noindent
\textbf{Corner start} (first piece $= a$).

As established above, only $c$, $e$ and $f$ can be bad, and the only pair
that can be simultaneously bad is $\{c,e\}$; in particular, whenever $f$ is
bad, neither $c$ nor $e$ can be bad. The individual counts follow from a
symmetry argument: piece $c$ is bad precisely when it appears first among
$\{b,c,f\}$, and since each of these three pieces is equally likely to
appear first, this happens in $5!/3=40$ orderings. Likewise, $e$ is bad
precisely when it appears first among $\{b,d,e,f\}$, in $5!/4=30$
orderings, and $f$ is bad precisely when it appears first among
$\{c,e,f\}$, in $5!/3=40$ orderings. Subtracting the $14$ orderings of
Table~\ref{tab:corner_count} in which $c$ and $e$ are simultaneously bad
yields the counts in Table~\ref{tab:corner_onebad}.

\begin{table}[h]
\centering
\caption{Counting orderings in the corner-start case with exactly one bad piece.}
\label{tab:corner_onebad}
\renewcommand{\arraystretch}{1.2}
\begin{tabular}{c p{6.5cm} c}
\toprule
Bad piece & Condition & Count \\
\midrule
$c$ only &
$c\prec\{b,f\}$, excluding the $14$ orderings with $e$ also bad &
$40-14=26$ \\
$e$ only &
$e\prec\{b,d,f\}$, excluding the $14$ orderings with $c$ also bad &
$30-14=16$ \\
$f$ only &
$f\prec\{c,e\}$; here neither $c$ nor $e$ can be bad &
$40$ \\
\bottomrule
\end{tabular}
\end{table}

Total: $26+16+40=82$ orderings. One corner piece contributes $82$
sequences; four corner pieces give $4\times 82=328$.

\medskip\noindent
\textbf{Middle start} (first piece $= b$).

Only $d$ and $f$ can be bad. By the same symmetry argument, $d$ is bad
precisely when it appears first among $\{a,d,e\}$, in $5!/3=40$ orderings,
and $f$ is bad precisely when it appears first among $\{c,e,f\}$, in
$5!/3=40$ orderings. In $16$ orderings both are bad, as computed in Table \ref{tab:middle_start}.
Hence the number of orderings with exactly one bad piece is
\[
(40-16)+(40-16)=24+24=48.
\]
One middle piece contributes $48$ sequences; two middle pieces give
$2\times 48=96$. Combining both symmetry classes gives a total of
$328+96=424$ assembly sequences with exactly one additional disconnected
component.

\section{Counting assembly sequences with zero disconnected components}

Now we count the assembly sequences in which no bad piece occurs, so that the partial assembly remains connected throughout. Equivalently, every piece placed
after the seed must be adjacent to at least one piece that has already been placed.  A useful observation shortens the count: in every branch below, once four pieces have been placed, each of the two remaining pieces is adjacent to a placed piece, so the sequence can always be completed in exactly $2!=2$ ways.

\medskip\noindent
\textbf{Corner start} (first piece $= a$, neighbours $b,d$).

The second piece must be $b$ or $d$. Table~\ref{tab:corner_connected}
enumerates the possibilities for the third piece in each case, together
with the admissible fourth pieces.

\begin{table}[h]
\centering
\caption{Counting connected orderings in the corner-start case by the second and third pieces.}
\label{tab:corner_connected}
\renewcommand{\arraystretch}{1.2}
\begin{tabular}{c c p{5.5cm} c}
\toprule
Second piece & Third piece & Admissible fourth pieces & Count \\
\midrule
$b$ & $c$ & $d,e,f$: $3$ choices, then $2$ completions & $3\times2=6$ \\
$b$ & $d$ & $c,e$: $2$ choices, then $2$ completions & $2\times2=4$ \\
$b$ & $e$ & $c,d,f$: $3$ choices, then $2$ completions & $3\times2=6$ \\
$d$ & $b$ & $c,e$: $2$ choices, then $2$ completions & $2\times2=4$ \\
$d$ & $e$ & $b,f$: $2$ choices, then $2$ completions & $2\times2=4$ \\
\bottomrule
\end{tabular}
\end{table}

Total: $6+4+6+4+4=24$ orderings. One corner piece contributes $24$
sequences; four corner pieces give $4\times 24=96$.

\medskip\noindent
\textbf{Middle start} (first piece $= b$, neighbours $a,c,e$).

The second piece must be $a$, $c$ or $e$. If the second piece is $a$, the
placed set $\{a,b\}$ coincides with that of the corner-start case with
second piece $b$, so there are $6+4+6=16$ completions; by the left--right
symmetry of the grid, the second piece $c$ likewise yields $16$. If the second piece is $e$, then any ordering of the remaining $4$ corner pieces are admissible since each corner piece is connected to one edge piece $b$ or $e$. The second piece $e$ thus yields $4\times6=24$ completions, and one middle piece contributes $16+16+24=56$ sequences in total; two middle pieces give
$2\times 56=112$. Combining both symmetry classes gives a total of
$96+112=208$ connected assembly sequences.

In summary, for the \(3\times2\) puzzle, the complete distribution of assembly sequences according to the number of additional disconnected components is
\(
6! = 720 = 208 + 424 + 88,
\)
where \(208\) assembly sequences remain connected throughout with zero disconnected components, \(424\) introduce one additional disconnected component, and \(88\) introduce two additional disconnected components.  




\end{document}